\documentclass[3p,12pt]{elsarticle}
\usepackage{amsthm}
\usepackage{amssymb}
\usepackage{amsmath}
\usepackage{epstopdf}
\usepackage{graphicx}
\usepackage{multirow}
\usepackage{subcaption}
\usepackage{listings}
\usepackage{array}
\usepackage{flafter}
\usepackage[section]{placeins}
\usepackage{epsfig}
\usepackage{booktabs}
\usepackage{longtable}
\usepackage{xcolor}
\usepackage[ruled]{algorithm2e}
\usepackage{color}
\usepackage{algpseudocode}
\numberwithin{equation}{section}

\newproof{pf}{Proof}[section]
\newdefinition{rmk}{Remark}
\definecolor{darkgreen}{rgb}{0,0.5,0}

\biboptions{sort&compress}

\journal{}
\begin{document}

\begin{frontmatter}

\title{FEALPy.OPT: A Multi-Backend Tensor Computing Framework for Metaheuristic Optimization on CPU and GPU Platforms}

\author[A]{Han-Bin Liu}
\author[A,B]{Huayi Wei\corref{cor1}}
\ead{weihuayi@xtu.edu.cn}
\cortext[cor1]{Corresponding author.}

\address[A]{School of Mathematics and Computational Science, Xiangtan University, Xiangtan 411105, China}
\address[B]{National Center of Applied Mathematics in Hunan, Hunan Key Laboratory for Computation
 and Simulation in Science and Engineering, Xiangtan 411105, China}

\begin{abstract}

Metaheuristic algorithms (MAs) are widely used to solve complex optimization problems, but their computational costs can become substantial as population size and problem dimensionality increase. This paper presents FEALPy.OPT, an open-source Python framework for metaheuristic optimization built on the FEALPy computational infrastructure. It integrates a unified optimizer abstraction with multi-backend tensor computation, supporting NumPy and PyTorch backends on CPU and GPU platforms. To evaluate its computational performance, five metaheuristic algorithms are tested on ten benchmark functions under three problem scales, with population sizes ranging from 100 to 10,000 and dimensions from 30 to 100. The results show that NumPy/CPU achieves the shortest average runtime for small-scale problems, whereas PyTorch/GPU becomes increasingly advantageous as the problem scale grows. For large-scale problems, PyTorch/GPU is approximately 9.50 times faster than PyTorch/CPU and 32.96 times faster than NumPy/CPU. These findings demonstrate the value of multi-backend and multi-platform support for improving the computational flexibility and efficiency of metaheuristic optimization.

\end{abstract}

\begin{keyword}
Metaheuristic optimization; Python framework; Tensor computing; Multi-backend; GPU acceleration

\end{keyword}

\end{frontmatter}

\section{Introduction}

In real-world scenarios, numerous constrained and unconstrained optimization problems necessitate resolution \cite{LingyunSanyang2023}.
Within the field of optimization, metaheuristic algorithms (MAs) have emerged as efficient and powerful tools for tackling complex issues \cite{LingyunSanyang2024}.
Over the past decade, MAs have garnered significant attention from researchers and engineers across various applications.
There are four primary reasons for the increasing popularity of these algorithms \cite{NirmalAliSanat2021,AmirPouyaKazemSeyedali2022}.
First, most MAs are inspired by simple natural concepts, making them easy to modify and adapt for specific applications.
The second reason is their black-box property; they do not require derivative information.
Unlike gradient-based algorithms, MAs do not necessitate knowledge of the gradient or Hessian of the objective function, which can often be difficult or impossible to ascertain.
Thirdly, most MAs can accommodate discrete variables, whereas many theoretical optimization methods struggle due to their reliance on the derivatives of the objective function.
Finally, and perhaps most importantly, MAs possess a degree of randomness that enables them to escape local optima stagnation.

Despite these advantages, MAs may incur substantial computational costs when applied to high-dimensional optimization problems or when large populations of candidate solutions are required.
This challenge is particularly relevant to population-based MAs, which constitute a broad family encompassing evolutionary algorithms, swarm-based methods, and other optimization approaches. particle swarm optimization (PSO) \cite{YansongYanminXiaoyan2025}, Genetic algorithms (GAs) \cite{AijiaYinshengYanmin2021}, and differential evolution (DE) \cite{XuLiAijia2018} are representative examples of this family.
In recent years, numerous additional population-based algorithms have been proposed, including the Snow Ablation Optimizer \cite{LingyunSanyang2023A}, Catch Fish Optimization Algorithm \cite{HemingQixianYuhao2024}, Crayfish Optimization Algorithm \cite{HemingHonghuaChangsheng2023}, and Genghis Khan shark optimizer \cite{GangYuxuanGuo2023}.
Although these algorithms employ different search mechanisms, many share a common computational pattern in which candidate solutions are repeatedly generated, evaluated, and updated throughout the optimization process.
As the population size, problem dimensionality, or computational cost of objective function evaluation increases, the overall workload can become substantial.
Consequently, computational efficiency depends not only on the underlying search strategy but also on the implementation of the associated numerical operations.
In particular, implementations that rely heavily on explicit loops over candidate solutions or decision variables may introduce additional overhead and underutilize modern computing hardware.
These considerations motivate the development of efficient computational frameworks that organize optimization operations effectively and facilitate the utilization of diverse computing resources.

Conventional implementations of MAs often rely on explicit loops to update individual search agents and their decision variables.
The positions of the search agents are typically organized into an \(n \times d\) matrix, where \(n\) denotes the population size and \(d\) represents the dimensionality of the search space.
In a loop-based implementation, the elements of this matrix are updated sequentially based on the current population and randomly sampled values.
Such implementations may introduce considerable overhead when the population size or problem dimensionality increases.
Vectorization provides an effective way to alleviate this overhead by replacing element-wise loops with batch array operations.
In Python, NumPy provides efficient array operations that facilitate the vectorized implementation of population-based MAs \cite{AmirPouyaKazemSeyedali2022}.
However, NumPy \cite{CharlesJarrodStefan2020} primarily supports CPU-based computation, whereas PyTorch \cite{Paszke2019PyTorch} provides a tensor computing interface that supports both CPU and GPU execution.
By exploiting the massive parallel computing capabilities of GPUs, tensor-based implementations may achieve substantial performance improvements for sufficiently large computational workloads, although the actual benefits depend on the problem characteristics and execution overhead.
Therefore, vectorization alone does not fully address the computational requirements of large-scale metaheuristic optimization.
Building on vectorized algorithm implementations, a computational framework that supports multiple backends and computing platforms can provide greater flexibility in utilizing available hardware resources and improving computational efficiency.

Existing Python libraries have facilitated the implementation and application of MAs by providing collections of optimization methods and reusable computational components \cite{AmirPouyaKazemSeyedali2022,evolopy,deap,pyswram,mealpy,niapy}.
Some frameworks focus on integrating a broad range of MAs, while others emphasize vectorized implementations to improve computational efficiency.
Meanwhile, tensor computing libraries provide mechanisms for executing numerical operations on different computing devices.
However, these capabilities address different aspects of metaheuristic optimization, and their integration into a unified framework remains an important consideration.
In particular, a framework that combines a consistent interface for MAs with multi-backend tensor computation and support for both CPU and GPU execution can provide a flexible foundation for algorithm development and computational performance evaluation.

To address these challenges, this paper presents FEALPy.OPT, an open-source Python framework for metaheuristic optimization that integrates algorithm implementation with multi-backend tensor computation and multi-platform execution.
Built on FEALPy's computational infrastructure, the framework aims to provide a unified and extensible environment for implementing, executing, and evaluating MAs \cite{fealpy}.
By supporting different numerical backends and computing devices, FEALPy.OPT facilitates the adaptation of optimization algorithms to diverse computational environments while maintaining a consistent programming interface.
The framework also provides a basis for systematically investigating the computational performance of metaheuristic algorithms under different backend and device configurations.
The main contributions of this work are summarized as follows:
\begin{itemize}
\item Framework design and implementation. We develop FEALPy.OPT, an open-source Python framework for metaheuristic optimization, providing an organized and extensible environment for implementing and using metaheuristic algorithms.

\item Multi-backend and multi-platform computing. We integrate backend-aware tensor operations into the optimization framework, enabling supported algorithms to execute with different numerical backends and computing devices, including CPU and GPU configurations.

\item Experimental evaluation. We conduct benchmark experiments to evaluate the computational performance of FEALPy.OPT under different backend and device configurations, examining the effects of these configurations on optimization runtime and computational efficiency.
\end{itemize}
This paper is structured as follows:
Section \ref{sec:related-work} gives a brief overview of the literature and similar works from researchers trying to make MAs faster and more efficient.
Section \ref{sec:fealpy.opt} provides a general framework of FEALPy.OPT.
Section \ref{sec:experiment} provides the experimental results and discussions.
Finally, the conclusion and future work is discussed in section \ref{sec:conclusion}

\section{Related work}\label{sec:related-work}

Research on MAs has primarily focused on improving their optimization capabilities, particularly convergence speed and global search performance.
Beyond algorithmic improvements, computational efficiency can also be enhanced through more efficient implementations and better utilization of hardware resources.
Although implementation-level optimization has received comparatively less attention in the metaheuristic literature, it is an important consideration, particularly for computationally expensive optimization problems.
The execution time of an algorithm depends not only on its search mechanism but also on how its computational operations are implemented and executed.
Therefore, several studies have investigated implementation strategies and computing techniques to accelerate metaheuristic optimization.

One important approach is parallel and distributed computing, which exploits multiple processing cores or computing nodes to execute independent tasks simultaneously \cite{Teodor2019}.
By distributing computational workloads across available resources, these approaches can reduce the overall execution time of optimization algorithms.
Parallel and distributed computing techniques have been investigated in evolutionary algorithms and other metaheuristic methods \cite{AlbaTomassini2002,AbtinJonatan2009}.
In addition to these approaches, vectorization provides another effective way to improve computational efficiency.
By organizing population-level operations into batch array computations, vectorized implementations can reduce the overhead associated with explicit loops and improve numerical computation efficiency.
Hemmasian et al. \cite{AmirPouyaKazemSeyedali2022} introduced VecMetaPy, a vectorized framework that utilizes NumPy to accelerate the implementation of several MAs.

While vectorization improves computational efficiency through optimized array operations, the choice of numerical backend and computing device also plays an important role in determining the performance of optimization algorithms.
Different tensor computing frameworks provide different computational capabilities and hardware support.
For example, NumPy \cite{CharlesJarrodStefan2020} primarily provides CPU-based numerical computation, whereas PyTorch \cite{Paszke2019PyTorch} supports tensor operations on both CPUs and GPUs.
By exploiting GPU computing capabilities, tensor-based implementations may offer performance advantages for sufficiently large computational workloads, although the actual benefits depend on problem characteristics and execution overhead.
These considerations motivate the development of optimization frameworks that integrate a unified algorithm interface with multiple computational backends and heterogeneous computing platforms, providing greater flexibility in utilizing available hardware resources.

This work is built upon FEALPy \cite{fealpy}, a cross-platform intelligent numerical simulation engine.
By leveraging its tensor computing infrastructure and backend management mechanism, we develop FEALPy.OPT as an extensible framework for metaheuristic optimization, integrating a unified algorithm interface with multi-backend tensor computation and multi-platform execution.

\section{FEALPy.OPT Framework Design}\label{sec:fealpy.opt}

\subsection{Overall Architecture}
\label{subsec:overall_architecture}

FEALPy.OPT is an optimization module developed within the FEALPy framework. It combines a shared optimizer abstraction and algorithm-specific implementations with FEALPy's tensor computing infrastructure. Its architecture comprises three interacting parts: the optimizer abstraction, the algorithm layer, and the backend-managed tensor operations.

The optimizer abstraction is defined in \texttt{opt/optimizer\_base.py}. The \texttt{Optimizer} class provides common configuration and optimization-related utilities, while \texttt{opt\_alg\_options()} assembles problem and execution parameters into a dictionary. Algorithm implementations are organized into subpackages, including swarm-based, evolutionary-based, physics-based, bio-based, human-based, music-based, math-based, hybrid, and multi-objective methods. This organization supports code navigation and extension, but does not imply that algorithms within or across these groups share identical search procedures.

For numerical computation, the algorithms use FEALPy's backend manager, imported as \texttt{bm}, to access tensor operations. The manager dispatches these operations to the currently selected backend, separating algorithm code from backend-specific operation implementations. Actual compatibility still depends on the operations and dependencies used by each algorithm, as well as the selected execution environment.

\begin{figure}[h]
    \centering
    \includegraphics[width=0.96\textwidth]
        {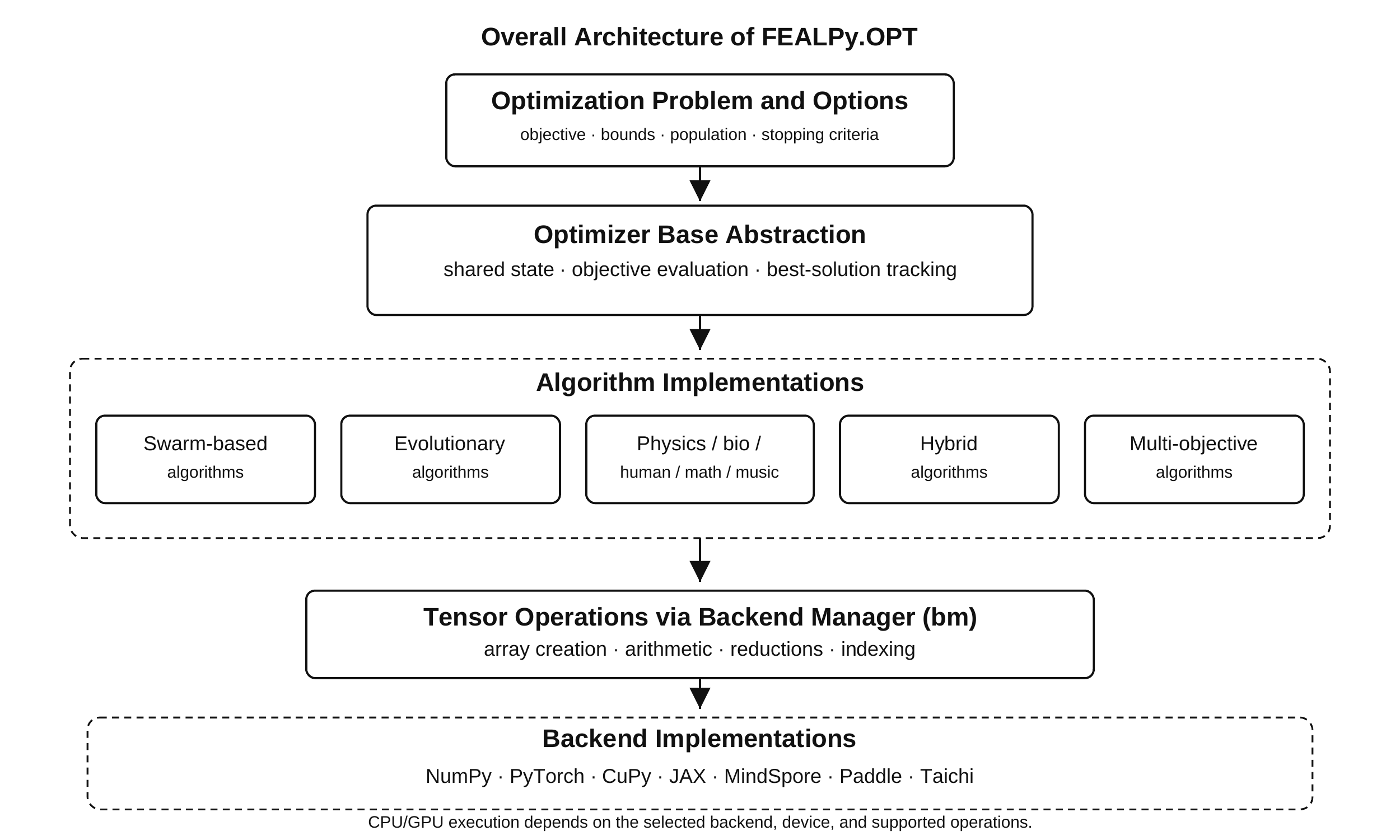}
    \caption{Overall architecture of FEALPy.OPT and its
    interaction with the FEALPy tensor computing infrastructure.}
    \label{fig:overall_architecture}
\end{figure}

Figure~\ref{fig:overall_architecture} illustrates the overall architecture of FEALPy.OPT and its interaction with the underlying FEALPy tensor computing infrastructure.

\subsection{Unified Optimizer Abstraction}
\label{subsec:unified_optimizer_abstraction}

FEALPy.OPT organizes its optimization algorithms around the \texttt{Optimizer} base class. The helper function \texttt{opt\_alg\_options()} assembles a configuration dictionary containing key information such as the initial solution or population \texttt{x0}, objective function \texttt{objective}, population size \texttt{NP}, dimension \texttt{ndim}, search domain \texttt{domain}, and iteration budget \texttt{MaxIters}. The dimension is derived from the last dimension of \texttt{x0}. During construction, the base class stores the configuration and initializes shared attributes, including the current solution representation, population size, dimension, optional domain bounds, and records used for convergence or population-diversity analysis.

The base class also provides reusable utilities, including \texttt{fun()} for invoking the configured objective function and maintaining an evaluation counter, and \texttt{update\_gbest()} for updating a global-best record. These utilities do not replace algorithm-specific state management: for example, individual algorithms initialize global-best information when required, and the presence of an option in the configuration does not guarantee that the base class enforces its associated stopping criterion.

The \texttt{run()} method in \texttt{Optimizer} raises \texttt{NotImplementedError}; subclasses define the actual search procedure, including initialization, state updates, and termination. Thus, PSO can maintain particle velocities and personal/global bests, whereas GA performs genetic operations and selection. The shared abstraction provides a common software interface and reusable utilities without imposing a uniform internal execution sequence on all algorithms.

\subsection{Tensor-based Population Computation}
\label{subsec:tensor_based_population_computation}

Population-based metaheuristic algorithms repeatedly update candidate solutions according to their search mechanisms. In representative continuous optimization algorithms in FEALPy.OPT, a population is typically organized as a two-dimensional tensor,
\begin{equation}
    \mathbf{X}
    =
    \begin{bmatrix}
        \mathbf{x}_1^\top \\
        \mathbf{x}_2^\top \\
        \vdots \\
        \mathbf{x}_N^\top
    \end{bmatrix}
    \in \mathbb{R}^{N \times D},
    \label{eq:population_tensor}
\end{equation}
where $N$ is the population size, $D$ is the decision-space dimension, and each row represents one candidate solution. Depending on the algorithm, additional tensors may store velocities, personal-best positions, or objective values. This is a typical representation in the inspected implementations, rather than a mandatory data structure enforced for every optimizer.

The PSO implementation illustrates population-level tensor computation: velocity and position updates use the current population, personal-best positions, the global-best solution, and random tensors; boundary handling is also performed through tensor operations. The GA implementation uses random tensors and Boolean masks for crossover and mutation, while QPSO includes population reductions and tensor-based position updates. These examples show that multiple algorithm families use tensor operations to process population-level numerical quantities, although their search procedures and control flows remain algorithm-specific.

The optimizer base class provides the \texttt{fun(x)} interface, which forwards its input to the configured objective function. This allows an algorithm to submit a population in a single call when the objective function accepts that input shape. It does not guarantee that every objective function internally evaluates candidates in a vectorized manner. Likewise, tensor-based representation should not be equated with complete vectorization, and vectorization alone does not establish GPU acceleration. The execution of these operations on different backends and devices depends on the underlying infrastructure, as discussed next.

\subsection{Backend and Device Abstraction}
\label{subsec:backend_device_abstraction}

\begin{figure}[h]
    \centering
    \includegraphics[width=\textwidth]
        {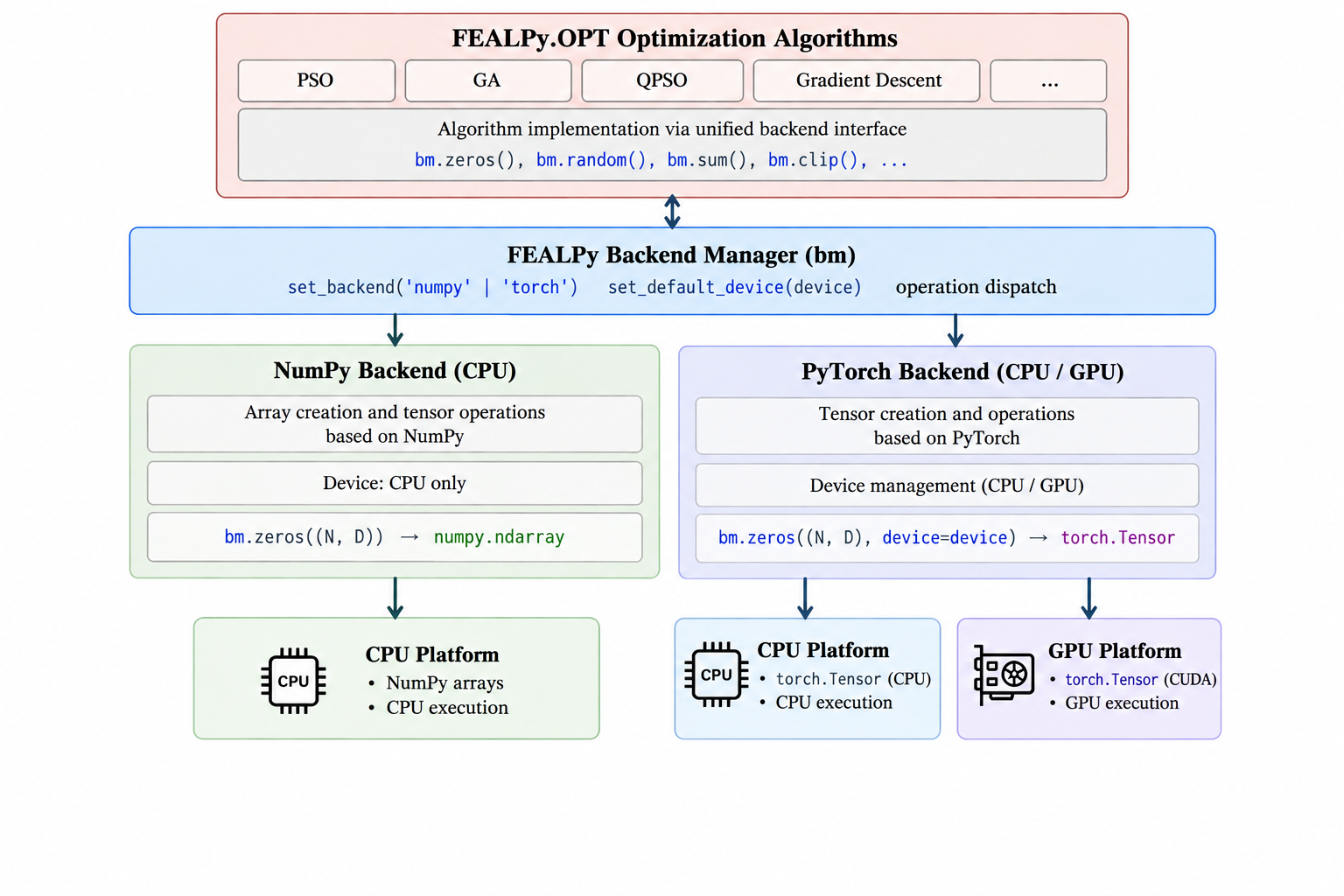}
    \caption{Multi-backend and multi-platform execution in
    FEALPy.OPT. Optimization algorithms invoke tensor operations
    through the backend manager (\texttt{bm}), which dispatches
    operations to the selected NumPy or PyTorch backend.
    The backends provide different device capabilities:
    NumPy supports CPU execution, whereas PyTorch provides
    interfaces for CPU and GPU execution.}
    \label{fig:backend_device_abstraction}
\end{figure}

FEALPy.OPT accesses tensor operations through FEALPy's backend manager, imported as \texttt{bm}. The manager loads a backend when required, selects the active backend through \texttt{set\_backend()}, and forwards attribute and method access to the corresponding backend proxy. Calls such as \texttt{bm.zeros()}, \texttt{bm.sum()}, and \texttt{bm.clip()} can therefore be expressed without directly invoking a specific tensor library at each call site.
Figure~\ref{fig:backend_device_abstraction} illustrates how this backend abstraction connects optimization algorithms to different tensor backends and execution platforms.

Backend selection changes the implementation used for subsequent \texttt{bm} operations; it does not automatically convert tensors that have already been created. For example, switching from NumPy to PyTorch does not turn an existing NumPy array into a PyTorch tensor. Input populations and algorithm-specific states must therefore be compatible with the selected backend, and backend switching should not be interpreted as automatic migration of an optimizer's data.

Device handling differs between backend implementations. The NumPy backend identifies its tensors as CPU-based, does not support setting a default device through \texttt{set\_default\_device()}, and restricts \texttt{device\_put()} to the CPU or an unspecified device. The PyTorch backend, by contrast, provides interfaces for setting a default device, querying a tensor's device, transferring tensors with \texttt{device\_put()}, and specifying a device in relevant tensor-creation functions. These mechanisms allow PyTorch operations to use supported devices, including GPUs when the installed PyTorch environment and hardware permit it. The existence of these interfaces does not, by itself, demonstrate successful GPU execution of every optimization algorithm.

Representative algorithms such as PSO, QPSO, and GA invoke tensor operations through \texttt{bm}, providing source-level evidence that their numerical computations use the backend abstraction. However, use of this interface is distinct from verified cross-backend compatibility. Correct execution under a particular configuration also depends on input data, internal states, objective-function behavior, and the operations required by the algorithm. Accordingly, the source code establishes the shared operation-dispatch mechanism and backend-specific device interfaces, whereas compatibility claims for individual algorithms require algorithm-level validation. This distinction provides the appropriate basis for interpreting the execution results reported in the experimental study.

\section{Experiment and analysis}\label{sec:experiment}

The experiments were conducted on a computer equipped with
an Intel Core i7-11800H processor, 16 GB of RAM, and an
NVIDIA GeForce RTX 3060 GPU. To evaluate the computational
performance of FEALPy.OPT across different computational
backends and platforms, three execution environments were
considered: NumPy on CPU, PyTorch on CPU, and PyTorch on GPU.

\begin{table}[!h]
\centering
\caption{Benchmark functions used in the experiments.}
\label{tab:benchmark-functions}
\begin{tabular}{llll}
  \hline
  Function & Name & Expression & Range \\ \hline
  {F1}&{Step}&\(F=\sum\limits_{i=1}^{D}(x_i+0.5)^2\)&{[-5.12,5.12]}\\
  {F2}&{Sphere}&\(F=\sum\limits_{i=1}^{D}x_i^2\)&{[-100,100]}\\
  {F3}&{Sum Squares}&\(F=\sum\limits_{i=1}^{D}i x_i^2\)&{[-10,10]}\\
  {F4}&{Bent Cigar}&\(F=x_1^2+10^6\sum\limits_{i=2}^{D}x_i^2\)&{[-100,100]}\\
  {F5}&{Schwefel 2.22}&\(F=\sum\limits_{i=1}^{D}|x_i|+\prod\limits_{i=1}^{D}|x_i|\)&{[-10,10]}\\
  {F6}&{Dixon--Price}&\(F=(x_1-1)^2+\sum\limits_{i=2}^{D}i(2x_i^2-x_i-1)^2\)&{[-10,10]}\\
  {F7}&{Rastrigin}&\(F=\sum\limits_{i=1}^{D}[x_i^2-10\cos(2\pi x_i)+10]\)&{[-5.12,5.12]}\\
  {F8}&{Rosenbrock}&\(F=\sum\limits_{i=1}^{D-1}[100(x_{i+1}-x_i^2)^2+(x_i-1)^2]\)&{[-30,30]}\\
  {F9}&{Griewank}&\(F=1+\frac{1}{4000}\sum\limits_{i=1}^{D}x_i^2-\prod\limits_{i=1}^{D}\cos\left(\frac{x_i}{\sqrt{i}}\right)\)&{[-600,600]}\\
  {F10}&{Ackley}&\(\begin{gathered}F=-20\exp\left(-0.2\sqrt{\frac{1}{D}\sum\limits_{i=1}^{D}x_i^2}\right)\\{}-\exp\left(\frac{1}{D}\sum\limits_{i=1}^{D}\cos(2\pi x_i)\right)+20+e\end{gathered}\)&{[-32,32]}\\
  \hline
\end{tabular}
\end{table}

Five metaheuristic optimization algorithms, namely the
Exponential Trigonometric Optimization Algorithm (ETOA) \cite{TranSamirMinh2024},
Differentiated Creative Search (DCS) \cite{PoominKhamronSirapat2024}, Particle Swarm
Optimization (PSO) \cite{YansongYanminXiaoyan2025}, Crayfish Optimization Algorithm (COA) \cite{HemingHonghuaChangsheng2023},
and Rime Optimization Algorithm (ROA) \cite{HangDongAli2023}, were selected for
the experiments. These algorithms were evaluated using the
ten benchmark functions presented in
Table~\ref{tab:benchmark-functions}. To investigate the
computational efficiency of the framework at different
problem scales, three combinations of population size and
problem dimension were considered: $(NP,D)=(100,30)$,
$(1000,50)$, and $(10000,100)$. The maximum number of
iterations was set to 1000 for all algorithms and
experimental configurations. The detailed experimental
settings are summarized in Table~\ref{tab:experimental-settings}.
Each configuration was independently executed 10 times
to obtain statistically representative runtime results.

\begin{table}[!h]
\centering
\caption{Experimental parameter settings.}
\label{tab:experimental-settings}
\begin{tabular}{ll}
\hline
Parameter & Setting \\ \hline
Optimization algorithms & ETOA, DCS, PSO, COA, ROA \\
Benchmark functions & F1--F10 \\
Population size and dimension
& (100, 30), (1000, 50), (10000, 100) \\
Maximum iterations & 1000 \\
Independent runs & 10 \\
Computational environments
& NumPy-CPU, PyTorch-CPU, PyTorch-GPU \\ \hline
\end{tabular}
\end{table}

The runtime of each algorithm was measured after optimizer
initialization, so that the reported time primarily reflected
the computational cost of the optimization process. For
CPU-based experiments, the runtime was measured using
\texttt{time.perf\_counter()}. For GPU-based experiments,
CUDA synchronization operations were performed immediately
before and after each optimization run to account for the
asynchronous execution of GPU operations and ensure accurate
end-to-end runtime measurements. The mean runtime, median
runtime, and standard deviation were calculated over the
10 independent runs. Furthermore, the GPU speedup ratio
relative to the PyTorch CPU implementation was calculated as
\begin{equation}
S_{\mathrm{GPU}} =
\frac{T_{\mathrm{PyTorch,CPU}}}
     {T_{\mathrm{PyTorch,GPU}}},
\end{equation}
where $T_{\mathrm{PyTorch,CPU}}$ and
$T_{\mathrm{PyTorch,GPU}}$ denote the mean runtimes of the
PyTorch implementations on CPU and GPU, respectively.
A speedup ratio greater than one indicates that GPU execution
is faster than CPU execution.

\begingroup
\scriptsize
\setlength{\tabcolsep}{2.5pt}
\setlength{\LTpre}{0pt}
\setlength{\LTpost}{0pt}
\renewcommand{\arraystretch}{0.85}
\begin{longtable}{l c c c c c c c}
    \caption{Mean runtime (s) of the five algorithms for each benchmark function, population size, dimension, and computational environment, averaged over 10 independent runs. The fastest environment for each algorithm and setting is shown in bold.}
    \label{tab:backend_performance}\\
    \toprule
    Function & $(NP,D)$ & Backend/device & ETOA & DCS & PSO & COA & ROA \\
    \midrule
    \endfirsthead
    \caption[]{Mean runtime (s) of the five algorithms for each benchmark function, population size, dimension, and computational environment, averaged over 10 independent runs. The fastest environment for each algorithm and setting is shown in bold (continued).}\\
    \toprule
    Function & $(NP,D)$ & Backend/device & ETOA & DCS & PSO & COA & ROA \\
    \midrule
    \endhead
    \midrule
    \multicolumn{8}{r}{Continued on next page}\\
    \endfoot
    \bottomrule
    \endlastfoot
    F1 & $(100,30)$ & NumPy/CPU & \textbf{0.4156346} & \textbf{0.3361387} & \textbf{0.1399642} & \textbf{0.3303344} & \textbf{0.151487} \\
       &            & PyTorch/CPU & 1.1141328 & 1.1491211 & 0.3347351 & 1.0100786 & 0.5916224 \\
       &            & PyTorch/GPU & 3.3517647 & 3.8302513 & 1.0595706 & 2.6902195 & 1.8040582 \\ \cline{2-8}
    \addlinespace[2pt]
       & $(1000,50)$ & NumPy/CPU & 3.7440569 & \textbf{2.0790158} & \textbf{0.9810107} & 3.2219774 & \textbf{1.1825369} \\
       &             & PyTorch/CPU & \textbf{3.1712396} & 3.0563929 & 1.0860317 & \textbf{2.4065308} & 1.2447647 \\
       &             & PyTorch/GPU & 3.2722612 & 3.7773109 & 1.0240027 & 2.6953906 & 1.7717118 \\ \cline{2-8}
    \addlinespace[2pt]
       & $(10000,100)$ & NumPy/CPU & 142.7291842 & 75.4850929 & 46.5610499 & 115.8561641 & 47.7185973 \\
       &               & PyTorch/CPU & 33.1833994 & 19.2972698 & 11.9610161 & 25.1353522 & 13.35639 \\
       &               & PyTorch/GPU & \textbf{3.5200863} & \textbf{3.6779402} & \textbf{1.1376501} & \textbf{2.7344361} & \textbf{1.8175377} \\
    \midrule
    F2 & $(100,30)$ & NumPy/CPU & \textbf{0.3514009} & \textbf{0.2975415} & \textbf{0.1127729} & \textbf{0.3145508} & \textbf{0.1492403} \\
       &            & PyTorch/CPU & 1.2390311 & 1.1192239 & 0.3203474 & 0.9627541 & 0.5726541 \\
       &            & PyTorch/GPU & 3.1505538 & 3.4647915 & 1.0176129 & 2.5600019 & 1.7117537 \\ \cline{2-8}
    \addlinespace[2pt]
       & $(1000,50)$ & NumPy/CPU & 3.3774757 & \textbf{1.8466224} & \textbf{0.9121242} & 2.9198484 & \textbf{0.9804193} \\
       &             & PyTorch/CPU & 3.2877105 & 3.3342429 & 1.1318922 & 2.6292714 & 1.5227332 \\
       &             & PyTorch/GPU & \textbf{3.2026274} & 3.6317131 & 1.0037661 & \textbf{2.5517178} & 1.780217 \\ \cline{2-8}
    \addlinespace[2pt]
       & $(10000,100)$ & NumPy/CPU & 125.4935131 & 73.0764459 & 44.6043904 & 114.4792837 & 45.803036 \\
       &               & PyTorch/CPU & 42.4315087 & 25.4168536 & 15.6154593 & 28.7043241 & 13.1516913 \\
       &               & PyTorch/GPU & \textbf{3.4541063} & \textbf{3.5828084} & \textbf{1.1363675} & \textbf{2.6128409} & \textbf{1.7487942} \\
    \midrule
    F3 & $(100,30)$ & NumPy/CPU & \textbf{0.3619363} & \textbf{0.2987334} & \textbf{0.1191394} & \textbf{0.316952} & \textbf{0.1495799} \\
       &            & PyTorch/CPU & 1.2558088 & 1.1696167 & 0.3425355 & 0.9931771 & 0.6044892 \\
       &            & PyTorch/GPU & 3.1291405 & 3.6753058 & 1.0502923 & 2.5698404 & 1.7508514 \\ \cline{2-8}
    \addlinespace[2pt]
       & $(1000,50)$ & NumPy/CPU & 3.3754545 & \textbf{1.8776071} & \textbf{0.9406971} & 2.931226 & \textbf{0.9908437} \\
       &             & PyTorch/CPU & 3.3291487 & 3.3305174 & 1.2035413 & 2.6615006 & 1.6292156 \\
       &             & PyTorch/GPU & \textbf{3.218235} & 3.7656908 & 1.0528519 & \textbf{2.5963433} & 1.8018737 \\ \cline{2-8}
    \addlinespace[2pt]
       & $(10000,100)$ & NumPy/CPU & 127.562729 & 75.1348701 & 46.7100134 & 116.3753405 & 47.974753 \\
       &               & PyTorch/CPU & 42.3895177 & 25.5498414 & 15.7043628 & 29.42303 & 13.2877055 \\
       &               & PyTorch/GPU & \textbf{3.6018291} & \textbf{3.7663126} & \textbf{1.1996569} & \textbf{2.6265211} & \textbf{1.8147987} \\
    \midrule
    F4 & $(100,30)$ & NumPy/CPU & \textbf{0.3547481} & \textbf{0.296689} & \textbf{0.1157622} & \textbf{0.315815} & \textbf{0.1485799} \\
       &            & PyTorch/CPU & 1.2906095 & 1.1970182 & 0.3543871 & 1.024254 & 0.5994463 \\
       &            & PyTorch/GPU & 3.1859549 & 3.5493655 & 1.0866858 & 2.6303403 & 1.8126082 \\ \cline{2-8}
    \addlinespace[2pt]
       & $(1000,50)$ & NumPy/CPU & 3.3774639 & \textbf{1.8918698} & \textbf{0.9218479} & 2.9438336 & \textbf{0.9891815} \\
       &             & PyTorch/CPU & 3.3646871 & 3.3998078 & 1.2059337 & 2.6494693 & 1.5947329 \\
       &             & PyTorch/GPU & \textbf{3.2494997} & 3.769118 & 1.1232949 & \textbf{2.649384} & 1.8262289 \\ \cline{2-8}
    \addlinespace[2pt]
       & $(10000,100)$ & NumPy/CPU & 125.732967 & 73.286634 & 44.9751808 & 114.9884333 & 46.0370881 \\
       &               & PyTorch/CPU & 43.3914172 & 25.3631778 & 15.4443748 & 29.8187452 & 13.6330342 \\
       &               & PyTorch/GPU & \textbf{3.608896} & \textbf{3.6834418} & \textbf{1.1425412} & \textbf{2.5780466} & \textbf{1.7896713} \\
    \midrule
    F5 & $(100,30)$ & NumPy/CPU & \textbf{0.3629413} & \textbf{0.3034847} & \textbf{0.1211367} & \textbf{0.3190321} & \textbf{0.149606} \\
       &            & PyTorch/CPU & 1.313116 & 1.1892116 & 0.3407049 & 0.9885971 & 0.5906532 \\
       &            & PyTorch/GPU & 3.0871929 & 3.5799491 & 1.0198846 & 2.5375471 & 1.7480494 \\ \cline{2-8}
    \addlinespace[2pt]
       & $(1000,50)$ & NumPy/CPU & 3.7117723 & \textbf{1.9739535} & \textbf{0.9609626} & 2.9930558 & \textbf{1.0116517} \\
       &             & PyTorch/CPU & 3.4331537 & 3.4018484 & 1.1679475 & 2.6590874 & 1.5796451 \\
       &             & PyTorch/GPU & \textbf{3.1631124} & 3.7352689 & 1.0205516 & \textbf{2.5643631} & 1.7690029 \\ \cline{2-8}
    \addlinespace[2pt]
       & $(10000,100)$ & NumPy/CPU & 126.9409219 & 74.2297143 & 45.075921 & 115.4579659 & 46.0337852 \\
       &               & PyTorch/CPU & 42.5903962 & 26.1210514 & 16.303542 & 29.1147054 & 13.2008103 \\
       &               & PyTorch/GPU & \textbf{3.5453042} & \textbf{3.7568254} & \textbf{1.1485506} & \textbf{2.6633089} & \textbf{1.7022169} \\
    \midrule
    F6 & $(100,30)$ & NumPy/CPU & \textbf{0.3634489} & \textbf{0.3124559} & \textbf{0.1300479} & \textbf{0.3301455} & \textbf{0.1635689} \\
       &            & PyTorch/CPU & 1.4877663 & 1.2247818 & 0.4064882 & 1.0847049 & 0.6606518 \\
       &            & PyTorch/GPU & 3.3094171 & 3.7587506 & 1.2617312 & 2.8271586 & 2.049463 \\ \cline{2-8}
    \addlinespace[2pt]
       & $(1000,50)$ & NumPy/CPU & \textbf{3.4044775} & \textbf{1.9018601} & \textbf{0.9950972} & 2.9141316 & \textbf{1.0485793} \\
       &             & PyTorch/CPU & 3.6360247 & 3.5151168 & 1.3687538 & \textbf{2.825181} & 1.7676594 \\
       &             & PyTorch/GPU & 3.4390934 & 3.9402637 & 1.2899541 & 2.8703943 & 1.9854397 \\ \cline{2-8}
    \addlinespace[2pt]
       & $(10000,100)$ & NumPy/CPU & 134.7232787 & 82.1190844 & 54.0531449 & 123.0816743 & 55.0908247 \\
       &               & PyTorch/CPU & 43.9303567 & 25.9680465 & 16.5531682 & 30.3800897 & 13.7977189 \\
       &               & PyTorch/GPU & \textbf{3.7039497} & \textbf{3.9188495} & \textbf{1.3284145} & \textbf{2.837836} & \textbf{1.9471667} \\
    \midrule
    F7 & $(100,30)$ & NumPy/CPU & \textbf{0.3814294} & \textbf{0.3354443} & \textbf{0.1406703} & \textbf{0.3303248} & \textbf{0.1732053} \\
       &            & PyTorch/CPU & 1.338952 & 1.5154535 & 0.4607674 & 1.1546293 & 0.7625454 \\
       &            & PyTorch/GPU & 3.1622711 & 3.6098507 & 1.1052204 & 2.6376306 & 1.7712839 \\ \cline{2-8}
    \addlinespace[2pt]
       & $(1000,50)$ & NumPy/CPU & 3.6920404 & \textbf{2.3001766} & 1.3296143 & 2.9891692 & \textbf{1.28767} \\
       &             & PyTorch/CPU & 3.5186109 & 3.4084606 & 1.2851704 & 2.7668466 & 1.6859494 \\
       &             & PyTorch/GPU & \textbf{3.2172736} & 3.7734436 & \textbf{1.1027514} & \textbf{2.6795655} & 1.8336013 \\ \cline{2-8}
    \addlinespace[2pt]
       & $(10000,100)$ & NumPy/CPU & 138.932928 & 88.6519493 & 61.1336935 & 123.791764 & 60.4710512 \\
       &               & PyTorch/CPU & 43.8014128 & 26.5185889 & 16.655472 & 30.1909244 & 13.8695001 \\
       &               & PyTorch/GPU & \textbf{3.6631022} & \textbf{3.7352479} & \textbf{1.2936815} & \textbf{2.7421717} & \textbf{1.9521193} \\
    \midrule
    F8 & $(100,30)$ & NumPy/CPU & \textbf{0.3744512} & \textbf{0.3199279} & \textbf{0.1294725} & \textbf{0.3207089} & \textbf{0.1599032} \\
       &            & PyTorch/CPU & 1.3739096 & 1.192486 & 0.3822181 & 1.0511143 & 0.6346988 \\
       &            & PyTorch/GPU & 3.2350256 & 3.7697709 & 1.2612256 & 2.884684 & 1.9977721 \\ \cline{2-8}
    \addlinespace[2pt]
       & $(1000,50)$ & NumPy/CPU & 3.5477365 & \textbf{1.8977491} & \textbf{1.0082223} & 2.9401548 & \textbf{1.0598342} \\
       &             & PyTorch/CPU & 3.6354436 & 3.504575 & 1.3467975 & 2.8204938 & 1.7207574 \\
       &             & PyTorch/GPU & \textbf{3.3474654} & 3.902367 & 1.2485032 & \textbf{2.7847386} & 1.9315945 \\ \cline{2-8}
    \addlinespace[2pt]
       & $(10000,100)$ & NumPy/CPU & 136.6276667 & 84.4101754 & 60.8553795 & 134.9183009 & 60.1571467 \\
       &               & PyTorch/CPU & 43.760019 & 26.4268398 & 16.691621 & 30.6668936 & 14.1930878 \\
       &               & PyTorch/GPU & \textbf{3.7479058} & \textbf{3.9595885} & \textbf{1.3396779} & \textbf{2.9366372} & \textbf{2.007642} \\
    \midrule
    F9 & $(100,30)$ & NumPy/CPU & \textbf{0.4000549} & \textbf{0.3276257} & \textbf{0.1654387} & \textbf{0.358464} & \textbf{0.1973706} \\
       &            & PyTorch/CPU & 1.4200746 & 1.6083027 & 0.5134383 & 1.2153165 & 0.8494504 \\
       &            & PyTorch/GPU & 3.4284379 & 3.7856786 & 1.210295 & 2.7553589 & 1.9622734 \\ \cline{2-8}
    \addlinespace[2pt]
       & $(1000,50)$ & NumPy/CPU & 3.9142002 & \textbf{2.2643551} & 1.8746106 & 3.6044631 & \textbf{1.5010061} \\
       &             & PyTorch/CPU & 3.6591911 & 3.4923277 & \textbf{1.3332665} & 2.8119183 & 1.7022304 \\
       &             & PyTorch/GPU & \textbf{3.3503147} & 3.8448635 & 1.5450089 & \textbf{2.7745213} & 1.9273161 \\ \cline{2-8}
    \addlinespace[2pt]
       & $(10000,100)$ & NumPy/CPU & 140.0127234 & 83.6653796 & 59.5901198 & 121.3916056 & 59.0248463 \\
       &               & PyTorch/CPU & 42.8120528 & 25.3416816 & 16.3759675 & 30.7158389 & 13.8979731 \\
       &               & PyTorch/GPU & \textbf{3.6451625} & \textbf{3.9880706} & \textbf{1.270796} & \textbf{2.8284732} & \textbf{1.9864022} \\
    \midrule
    F10 & $(100,30)$ & NumPy/CPU & \textbf{0.4055444} & \textbf{0.3263944} & \textbf{0.1594177} & \textbf{0.34017} & \textbf{0.1890693} \\
        &            & PyTorch/CPU & 1.4609257 & 1.6204518 & 0.585401 & 1.3411 & 0.907535 \\
        &            & PyTorch/GPU & 3.5435241 & 4.1259206 & 1.3517028 & 2.9680693 & 2.0603409 \\ \cline{2-8}
    \addlinespace[2pt]
        & $(1000,50)$ & NumPy/CPU & 3.8203387 & \textbf{2.129279} & 1.4130645 & 3.0897506 & \textbf{1.4220707} \\
        &             & PyTorch/CPU & 3.723792 & 3.5055 & 1.4406978 & \textbf{2.8066468} & 1.7393898 \\
        &             & PyTorch/GPU & \textbf{3.4567293} & 4.059222 & \textbf{1.3875694} & 2.8801181 & 2.0827565 \\ \cline{2-8}
    \addlinespace[2pt]
        & $(10000,100)$ & NumPy/CPU & 138.6424486 & 82.8435612 & 58.2056898 & 123.8878335 & 59.2963249 \\
        &               & PyTorch/CPU & 45.1463869 & 25.5111061 & 16.3578504 & 30.7197823 & 13.6849434 \\
        &               & PyTorch/GPU & \textbf{3.6433854} & \textbf{4.0427441} & \textbf{1.4009098} & \textbf{2.91} & \textbf{2.1278927} \\
\end{longtable}
\endgroup

Each entry in Table~\ref{tab:backend_performance} is the mean runtime over 10 independent runs for a fixed algorithm, benchmark function, $(NP,D)$ setting, and computational environment. To compare the environments at each scale, an additional arithmetic mean was calculated across the 50 algorithm--function entries (five algorithms and ten functions). For $(NP,D)=(100,30)$, these summary means for NumPy/CPU, PyTorch/CPU, and PyTorch/GPU were 0.263~s, 0.958~s, and 2.538~s, respectively. NumPy/CPU achieved the shortest runtime, being approximately 3.64 and 9.64 times faster than PyTorch/CPU and PyTorch/GPU, respectively. This indicates that, for relatively small workloads, the computational overhead associated with PyTorch and GPU execution may outweigh the benefits of their underlying computational capabilities. When the problem scale increased to $(NP,D)=(1000,50)$, the corresponding summary means were 2.190, 2.490, and 2.573~s, respectively. The performance differences among the three environments became considerably smaller, with NumPy/CPU maintaining a slight advantage. These results suggest that the benefits of GPU execution are not yet sufficient to compensate for its associated overhead at these two smaller problem scales.

A different performance pattern was observed for the largest problem scale, $(NP,D)=(10000,100)$. The summary means of the table entries for NumPy/CPU, PyTorch/CPU, and PyTorch/GPU were 87.678~s, 25.271~s, and 2.660~s, respectively. In particular, PyTorch/GPU was approximately 9.50 times faster than PyTorch/CPU and 32.96 times faster than NumPy/CPU. The comparison between the two PyTorch configurations highlights the benefit of GPU execution under a larger computational workload, while the comparison with NumPy/CPU further demonstrates the performance difference between the available computational environments. These results indicate that the relative efficiency of different backends and computing platforms varies with problem scale, and that GPU execution becomes particularly advantageous when the computational workload is sufficiently large.

\begin{figure}[h]
    \centering
    \includegraphics[width=\textwidth]
        {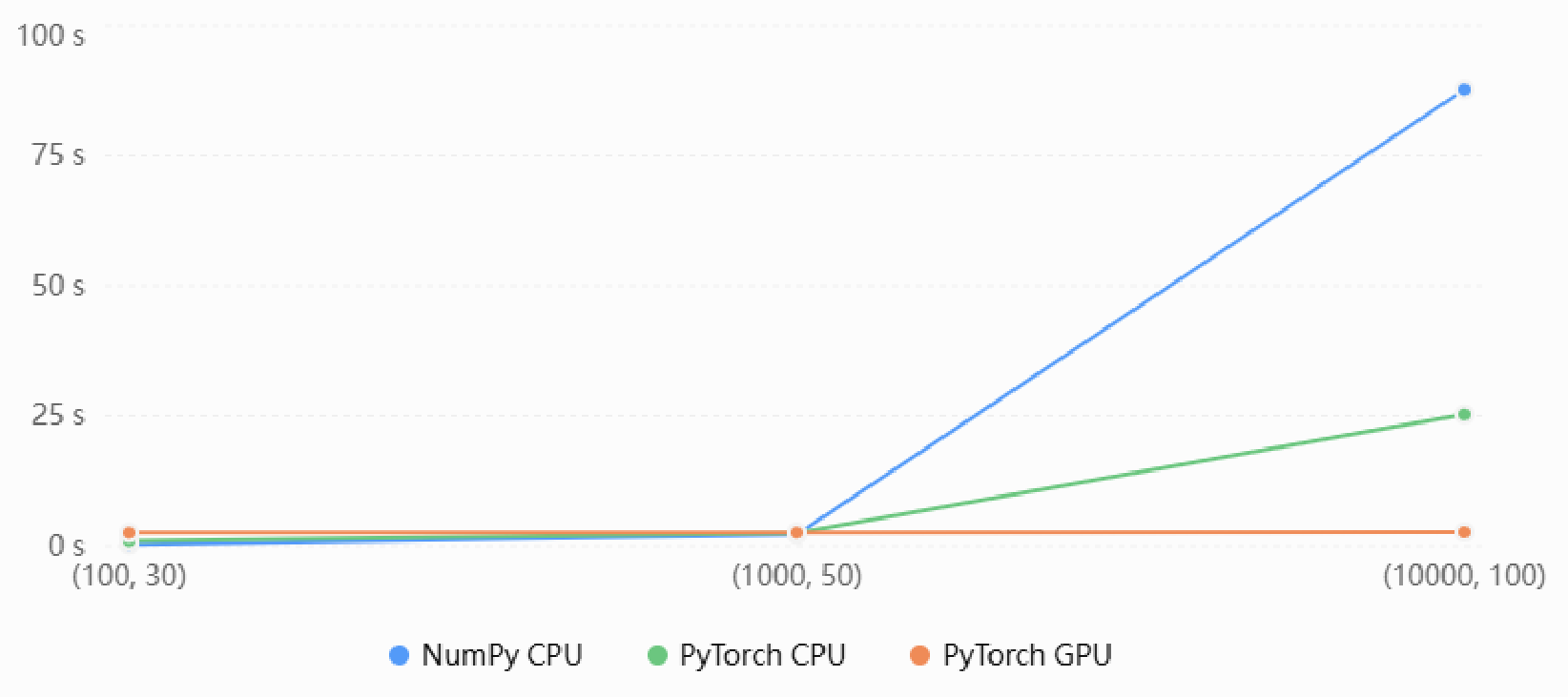}
    \caption{Comparison of average runtime across different computational environments and problem scales.}
    \label{fig:cpu_gpu_compare}
\end{figure}

To provide an overall comparison of the three computational environments, the mean runtime was further aggregated across the 50 algorithm–function combinations for each problem scale. Figure \ref{fig:cpu_gpu_compare} illustrates the variation in average runtime with increasing problem size. A clear difference can be observed in the performance trends of the three environments. While NumPy/CPU exhibits a substantial increase in runtime as the problem scale grows, PyTorch/GPU maintains a relatively stable runtime across the tested configurations. PyTorch/CPU shows an intermediate trend, with its runtime increasing at a slower rate than that of NumPy/CPU. These results indicate that the relative computational efficiency of different environments is closely related to the scale of the optimization workload.

Overall, the experimental results demonstrate the practical value of supporting multiple computational backends and platforms in FEALPy.OPT. Rather than relying on a single fixed computational environment, users can select a suitable configuration according to the scale of the optimization task. NumPy/CPU provides better performance for the smaller workloads considered in this study, whereas PyTorch/GPU offers substantial acceleration for the largest workload. PyTorch/CPU also provides an alternative execution environment, with performance between the other two configurations in the large-scale experiments. Therefore, the multi-backend and multi-platform design of FEALPy.OPT provides flexibility in utilizing available computational resources and enables more efficient execution across different problem scales.
\section{Conclusion and future work}\label{sec:conclusion}

This paper presents FEALPy.OPT, an open-source Python framework for metaheuristic optimization built on the FEALPy computational infrastructure. The framework integrates a unified optimizer abstraction with backend-aware tensor operations, providing a consistent interface for implementing and executing metaheuristic algorithms. By leveraging FEALPy's backend management mechanism, FEALPy.OPT supports NumPy and PyTorch backends, enabling execution on CPU and GPU platforms. This design separates algorithm implementations from backend-specific numerical operations and provides a flexible computational environment for metaheuristic optimization.

The computational performance of FEALPy.OPT was evaluated using five metaheuristic algorithms and ten benchmark functions under three problem scales. The experimental results demonstrate that the relative efficiency of different computational environments depends on the problem scale. NumPy/CPU performs best on small-scale problems, while the performance differences among the three environments become smaller at the medium scale. For large-scale problems, PyTorch/GPU achieves approximately 9.50 times and 32.96 times shorter average runtimes than PyTorch/CPU and NumPy/CPU, respectively. These results highlight the importance of selecting appropriate computational backends and devices according to the characteristics of optimization tasks, and demonstrate the practical value of integrating multi-backend tensor computation and multi-platform execution into a metaheuristic optimization framework.

Future work will focus on extending the algorithm coverage and improving backend compatibility across a broader range of metaheuristic methods. Further investigations will also consider additional tensor computing backends and computing platforms, as well as more diverse and computationally intensive optimization problems. These developments will help enhance the flexibility and applicability of FEALPy.OPT and provide further insights into the computational performance of metaheuristic optimization across heterogeneous computing environments.

{\bf CRediT authorship contribution statement}

{\bf Han-Bin Liu:} Methodology, Writing - original draft, Writing - review \& editing. {\bf Huayi Wei:} Software, Writing - original draft, Writing - review \& editing.
\\

{\bf Acknowledments}

This work is supported by the National Natural Science Foundation of China (Grant Nos. 12371410, 12261131501), the construction of innovative provinces in Hunan Province (Grant No. 2021GK1010), , the Postgraduate Scientific Research Innovation Project of Hunan Province (Grant No. CX20250944) and the Postgraduate Scientific Research Innovation Project of Xiangtan University (Grant No. XDCX2025Y205).
\\

{\bf Declaration of competing interest}

The authors declare that they have no known competing financial interests or personal relationships that could have appeared to influence the work reported in this paper.
\\

{\bf Data availability}

Data will be made available on request.

\end{document}